%% file: def-Tutte.tex
\documentclass[10pt,a4paper,reqno,english]{amsart} 
\usepackage{ae,aecompl}
\usepackage[cm]{aeguill}

\usepackage{pstricks,amsmath,here,latexsym,amssymb}
\usepackage{epsf,psfrag,epsfig,color,graphicx}
\usepackage[latin1]{inputenc}
\usepackage{t1enc}
\usepackage{fancyheadings} 
\usepackage{epic} 
\usepackage{eepic}

\newtheorem{prop}{Proposition} 
\newtheorem{lemma}[prop]{Lemma}

\newtheorem{thm}[prop]{Theorem}
\newcommand{\dem}{\noindent \textbf{Proof: }}

\theoremstyle{definition} 
\newtheorem{Def}[prop]{Definition}

\newcommand{\Ref}[1]{(\ref{#1})}
\newcommand{\findem}{\vspace{-.4cm} \begin{flushright} $\square~$ \end{flushright} \vspace{.4cm} }

\newcommand{\mG}{\mathcal{G}}
\newcommand{\mI}{\mathcal{I}}
\newcommand{\mE}{\mathcal{E}}
\newcommand{\mT}{\mathcal{T}}

\catcode`\@=11
\def\section{\@startsection{section}{1}%
 \z@{.7\linespacing\@plus\linespacing}{.5\linespacing}%
 {\normalfont\bfseries\scshape\centering}}

\def\subsection{\@startsection{subsection}{2}%
  \z@{.5\linespacing\@plus\linespacing}{.5\linespacing}%
  {\normalfont\bfseries\scshape}}

\def\subsubsection{\@startsection{subsubsection}{3}%
  \z@{.5\linespacing\@plus.7\linespacing}{-.5em}%
  {\normalfont\itshape}}
\catcode`\@=12
\newcommand{\ccirc}{}

\newcommand{\ite}{\noindent $\bullet~$}

\title[The Tutte polynomial for embedded graphs]
{A characterization of the Tutte polynomial via combinatorial embeddings}
\author{Olivier Bernardi}
\address{LaBRI, Universit\'e Bordeaux 1, 351 cours de la Lib\'eration, 33405 Talence Cedex, France}
\email{bernardi@labri.fr}
\date{\today}

\begin{document}
\maketitle

\begin{abstract}
We give a new characterization of the Tutte polynomial of graphs. Our
characterization is formally close (but inequivalent) to the original definition given by Tutte 
as the generating function of spanning trees counted according to \emph{activities}.
Tutte's notion of activity requires to choose a \emph{linear order} on
the edge set (though the generating function of the activities
 is, in fact, independent of this order). We define a new notion of
activity, the \emph{embedding-activity}, which requires to choose 
 a \emph{combinatorial embedding} of the graph, that is, a cyclic order of
the edges around each vertex. We prove 
that the Tutte polynomial  equals the generating function of spanning trees counted according to  \emph{embedding-activities} (this generating function being, in fact, independent of the embedding).\\ 
\end{abstract}

\section{Introduction}
In 1954, Tutte defined a graph invariant that he named
\emph{dichromate} because he thought of it as a bivariate
generalization of the \emph{chromatic polynomial}
\cite{Tutte:dichromate}. The first definition by Tutte was a generating function of spanning trees counted according to their \emph{activities} (see Theorem \ref{thm:activity-Tutte}). Since then, this polynomial, which is now known as the \emph{Tutte polynomial}, has been widely studied (see for instance, \cite{Brylawsky:Tutte-poly} and references therein).  We refer the reader to \cite[Chapter X]{Bollobas:Tutte-poly} for an easy-to-read but comprehensive survey of the properties and applications of the Tutte polynomial.\\


In this paper, we give a new characterization of the Tutte polynomial
of graphs. Our characterization is formally close to the original definition by
Tutte (compare Theorem \ref{thm:activity-Tutte} and Theorem
\ref{thm:activity-OB}) in terms of the \emph{activities} of spanning trees. 
Tutte's notion of activity  requires to choose a \emph{linear order} on the edge set. 
The Tutte polynomial is then the generating function of spanning trees counted according 
to their (internal and external) activities (this
generating function being, in fact, independent of the linear order;
see Theorem \ref{thm:activity-Tutte}). 
Our characterization of the Tutte polynomial requires instead to
choose an \emph{embedding} of the graph, that is, a cyclic order for the
incidences of edges around each vertex. Once the embedding is chosen, one
can define the (internal and external) \emph{embedding-activities} of
spanning trees. We prove that the Tutte polynomial is equal  to the
 generating function of spanning trees counted according 
to their (internal and external) embedding-activities 
(this generating function being, in fact, independent of the
embedding; see Theorem \ref{thm:activity-OB}).\\ 

Several other notions of \emph{activities} related to the Tutte polynomial have been introduced. In 1982 Las Vergnas gave another characterization of the Tutte polynomial as the generating function of orientations counted according to their \emph{cyclic-activities} \cite{Vergnas:Morphism-matroids-2}. Las Vergnas' notion of activity requires to linearly order the edge set. A connection with Tutte's notion of activity was established in \cite{Gioan-bij-tree-orientation}. An alternative notion of \emph{external activity} has also been introduced by Gessel and Wang \cite{Gessel-Wang:DFS} and further investigated in \cite{Gessel:Tutte-poly+DFS,Gessel:enumerative-csq-DFS}. This notion is related to the depth-first search algorithm and requires to linearly order the vertex set. \\

This paper is organized as follows. In Section
\ref{section:definition}, we introduce our notations about graphs, and define embeddings. In Section \ref{section:polyTutte}, we
recall some classical properties of the Tutte polynomial.  
In Section \ref{section:def-Tutte-OB}, we define the
\emph{embedding-activities}  and prove that the Tutte polynomial equals the generating function of the embedding-activities. 
Lastly, in Section \ref{section:remarks}, we mention some possible applications of our characterization of the Tutte polynomial, some of which will be developed in a forthcoming paper~\cite{OB:Tutte-plongement}.\\ 


\section{Combinatorial embeddings of graphs} \label{section:definition}
\begin{figure}[htb!]
\begin{center}
\input{exp-contraction.pstex_t}
\caption{Deletion and contraction of the edge $e=\{b,b'\}$.}\label{fig:exp-contraction}
\end{center}
\end{figure}
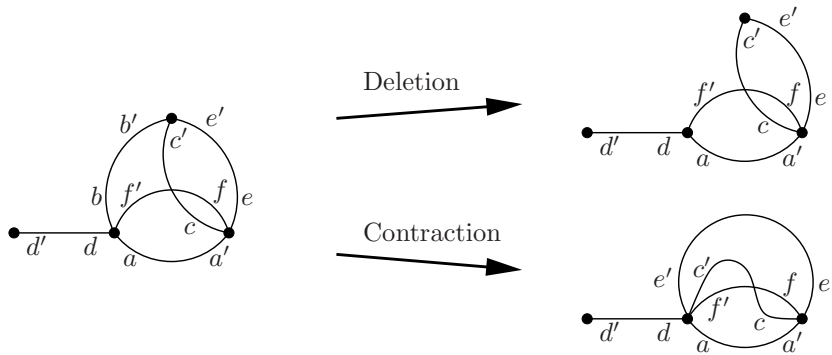

\noindent \textbf{Graphs. } In this paper we consider finite, undirected
graphs. Loops and multiple edges are allowed. Formally, a \emph{graph}
$G=(V,E)$ is a finite set of \emph{vertices} $V$, a finite set of
\emph{edges} $E$ and a relation of \emph{incidence} in $V\times E$
such that each edge $e$ is incident to either one  or two vertices. 
The \emph{endpoints} of an edge $e$ are the vertices incident to
$e$. A \emph{loop} is an edge with just one endpoint.  
If $e$ is an edge of $G$, we denote by $G_{\backslash e}$ the graph
obtained by \emph{deleting} $e$. If $e$ is not a loop, we denote by
$G_{/e}$ the  graph obtained by \emph{contracting} the edge $e$, that
is, by identifying its two endpoints (see Figure
\ref{fig:exp-contraction}). An \emph{isthmus} is an edge whose
deletion increases the number of connected components. \\

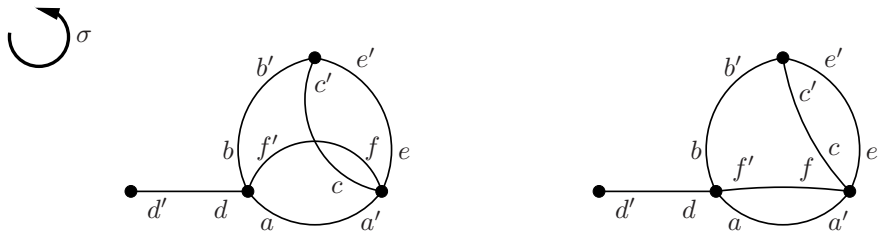
\begin{figure}[htb!]
\begin{center}
\input{exp-embedding.pstex_t}
\caption{Two embeddings of the same graph.}\label{fig:exp-embedding}
\end{center}
\end{figure}

\noindent \textbf{Embeddings. } 
We recall the notion of \emph{combinatorial map} introduced by Cori and Machi \cite{Cori:These-asterisque,Cori-Machi:survey}.  A \emph{combinatorial map} (or \emph{map} for short) $\mG=(H,\sigma,\alpha)$ is a set of \emph{half-edges} $H$, a permutation $\sigma$  and an involution without fixed point $\alpha$ on $H$ such that the group generated by $\sigma$ and $\alpha$ acts transitively on $H$. A map is \emph{rooted} if one of the half-edges is distinguished as the \emph{root}. For $h_0\in H$, we denote by $\mG=(H,\sigma,\alpha,h_0)$ the map $(H,\sigma,\alpha)$ rooted on $h_0$.
Two maps $\mG=(H,\sigma,\alpha)$ and $\mG'=(H',\sigma',\alpha')$ are \emph{isomorphic} if there is a bijection $\pi: H\mapsto H'$ such that $\pi\sigma=\sigma'\pi$ and $\pi\alpha=\alpha'\pi$. If $\mG$ and $\mG'$ are rooted on $h_0$ and $h_0'$ respectively, we also require that $\pi(h_0)=h_0'$.\\
Given a combinatorial map $\mG=(H,\sigma,\alpha)$, we consider the \emph{underlying} graph $G=(V,E)$, where $V$ is the set of cycles of $\sigma$, $E$ is the set cycles of $\alpha$ and the incidence relation is to have at least one common half-edge.  We represented the underlying graph of the map $\mG=(H,\sigma,\alpha)$ where the set of half-edges is $H=\{a,a',b,b',c,c',d,d',e,e',f,f'\}$, the involution $\alpha$  (written as a product of cycles of length 2) is $(a,a')$$(b,b')$$(c,c')$$(d,d')$$(e,e')$$(f,f')$  and the permutation $\sigma$ is $(a,f',b,d)$ $(d')$$(a',e,f,c)$$(e',b',c')$ on the left of Figure \ref{fig:exp-embedding}. 
Graphically, we represent the cycles of $\sigma$ by the counterclockwise order around each vertex. Hence, our drawing defines the map $\mG$ since the order around vertices gives the cycles of $\sigma$ (the involution $\alpha$ is immediately recovered from the edges). On the right of Figure \ref{fig:exp-embedding},  we represented the map $\mG'=(H,\sigma',\alpha)$, where $\sigma'=(a,f',b,d)$$(d')$$(a',e,c,f)$$(e',b',c')$. The maps  $\mG$ and $\mG'$ have isomorphic underlying graphs. \\
Note that the underlying graph of a map $\mG=(H,\sigma,\alpha)$ is always connected since  $\sigma$ and $\alpha$ act transitively on $H$. A (rooted) \emph{combinatorial embedding} (or \emph{embedding} for short) of a connected graph $G$ is a (rooted) combinatorial map $\mG$ whose underlying graph is isomorphic to $G$.  When an embedding $\mG$  of $G$ is given we shall write the edges of $G$  as pairs of half-edges (writing for instance $e=\{h,h'\}$). \\

Intuitively, a combinatorial embedding corresponds to the choice a cyclic order of the edges around each vertex. This order  can also be seen as a local planar embedding. In fact, combinatorial embeddings are closely related to the \emph{embeddings of  graphs on surfaces}. What we call \emph{surfaces} are compact 2-dimensional orientable surfaces without boundary (see \cite{Mohar:graphs-on-surfaces}). An \emph{embeddings of a graph $G$ in a surface $\mathcal{S}$} is a drawing of $G$ on $\mathcal{S}$ without intersecting edges such that each connected component of the complement of the graph (i.e. each \emph{face}) is homeomorphic to a disc. Consider the embeddings of $G$ in a surface $\mathcal{S}$ of genus $g$, defined up to orientation preserving homeomorphisms. An embedding on  $\mathcal{S}$ gives a cyclic order of the edges around each vertex, hence defines a combinatorial embedding. In fact, this relation is one-to-one. The embeddings of $G$ in a surface $\mathcal{S}$ of genus $g$ (defined up to orientation preserving homeomorphisms) are in one-to-one correspondence with the combinatorial embeddings $\mG=(H,\sigma,\alpha)$ of $G$ with \emph{Euler characteristic} $\chi_\mG=2-2g$ (up to isomorphisms) \cite[Thm. 3.2.4]{Mohar:graphs-on-surfaces}, where the Euler characteristic $\chi_\mG$ is defined as  the number of cycles of $\sigma$, plus the number of cycles of $\sigma \ccirc \alpha$, minus the number of cycles of $\alpha$.\\

In the following, we use the terms \emph{combinatorial maps} and \emph{embedded graphs} interchangeably. \emph{We do not require our embedded graphs to be planar} (i.e. to have an Euler characteristic equal to 2).\\

We now define the deletion and contraction of edges in embedded graphs.
Our definition  preserves the \emph{cyclic order} around each vertex. 
We represented the result of deleting and contracting the edge $e=\{b,b'\}$ in Figure \ref{fig:exp-contraction}. 

%
%
Let $G$ be a graph and  $\mathcal{G}=(H,\sigma,\alpha)$ one of its embeddings. Let $e=\{h_1,h_2\}$ be an edge of $G$. If $e$ is not an isthmus, we define the embeddings $\mG_{\backslash  e}$  of $G_{\backslash  e}$   by  
$\mG_{\backslash  e}=(H',\sigma',\alpha_{\backslash  e})$ where $H'=H\setminus\{h_1,h_2\}$, the involution $\alpha'$ is $\alpha$ restricted to $H'$ and 
\begin{eqnarray}\label{eq:deletion-plongement}
\begin{array}{l|ll}
\sigma_{\backslash e}(h)=
& \sigma\ccirc \sigma\ccirc \sigma(h)  & \textrm{if }  (\sigma(h)=h_1  \textrm{ and }\sigma(h_1)=h_2)  \textrm{ or } (\sigma(h)=h_2  \textrm{ and }\sigma(h_2)=h_1),\\
& \sigma\ccirc \sigma(h) & \textrm{if }  (\sigma(h)=h_1   \textrm{ and }\sigma(h_1)\neq h_2) \textrm{ or } (\sigma(h)=h_2  \textrm{ and }\sigma(h_2)\neq h_1),\\
& \sigma(h) & \textrm{otherwise}.\\
\end{array}
\end{eqnarray}
Similarly,  if $e$ is not a loop, we define the embeddings $\mG_{/ e}$  of $G_{/e}$  by $\mathcal{G}_{/e}=(H',\sigma',\alpha_{/e})$, where $H'=H\setminus\{h_1,h_2\}$, the involution $\alpha'$ is $\alpha$ restricted to $H'$ and 
\begin{eqnarray}\label{eq:contraction-plongement}
\begin{array}{l|ll}
\sigma_{/ e}(h)=
& \sigma\ccirc \sigma(h) & \textrm{if }  (\sigma(h)=h_1 \textrm{ and } \sigma(h_2)=h_2)  \textrm{ or }  (\sigma(h)=h_2 \textrm{ and } \sigma(h_1)=h_1),\\
& \sigma\ccirc \alpha\ccirc\sigma(h) & \textrm{if }  (\sigma(h)=h_1 \textrm{ and } \sigma(h_2)\neq h_2)  \textrm{ or }  (\sigma(h)=h_2 \textrm{ and } \sigma(h_1)\neq h_1), \\
&\sigma(h) & \textrm{otherwise}.\\
\end{array}
\end{eqnarray}


\vspace{.5cm}

\noindent \textbf{Spanning trees. }
Let $G=(V,E)$ be a graph. A spanning subgraph of $G$ is a graph  $G'=(V,E')$ where $E'\subseteq E$.
A spanning subgraph is entirely determined by its edge set and, by convenience, we will often identify the spanning subgraph with its edge set.  A \emph{tree} is a connected acyclic graph. A \emph{spanning tree} is a spanning subgraph which is a tree. \\ 
Let  $T$  be a spanning tree of $G$. An edge of $E$ is said to be \emph{internal} 
if it is in $T$ and \emph{external} otherwise. The
\emph{fundamental cycle} (resp. \emph{cocycle}) of an external
(resp. internal) edge $e$ is the set of edges $f$ such that the
spanning subgraph $T \setminus \{f\} \cup \{e\}$ (resp.  $T \setminus \{e\} \cup
\{f\}$) is a tree. Note that the fundamental cycle  (resp. cocycle) of
an external (resp. internal) edge $e$ contains only internal
(resp. external) edges apart from $e$. Moreover, if $e$ is internal
and $f$ external, then $e$ is in the fundamental cycle of $f$ if and
only if $f$ is in the fundamental cocycle of $e$. \\

Lastly, for any set $S$ we denote by $|S|$ its cardinality. 
If $S\subset S'$ and $s\in S'$, we write $S-s$ for $S\setminus\{s\}$ whether $s$ belongs to $S$ or not.

\section{The Tutte polynomial} \label{section:polyTutte}
\begin{Def} \label{def:Tutte-poly}
The \emph{Tutte polynomial} of a graph $G=(V,E)$ is
\begin{eqnarray}
T_G(x,y)=\sum_{S\subseteq E} (x-1)^{c(S)-c(G)}(y-1)^{c(S)+|S|-|V|}, \label{eq:Tutte-explicit}
\end{eqnarray}
where the sum is over all spanning subgraphs $S$ and $c(S)$ (resp. $c(G)$) denotes the number of connected components of $S$ (resp. of $G$).
\end{Def}

For example, if  $G$ is $K_3$ (the triangle) there are 8 spanning subgraphs.  The subgraph with no edge has contribution $(x-1)^2$, each subgraph with one edge has contribution  $(x-1)$, each subgraph with two edges has contribution  $1$ and the subgraph with three edges has contributions $(y-1)$. Summing up these contributions, we get $T_{K_3}(x,y)=(x-1)^2+3(x-1)+3+(y-1)=x^2+x+y$.\\


From Definition \ref{def:Tutte-poly}, it is easy to check that whenever $G$ is the disjoint union of two graphs $G=G_1 \cup G_2$, then 
$$T_G(x,y)~=~T_{G_1}(x,y)\times T_{G_2}(x,y).$$ 
This relation allows us to restrict our attention to connected graphs. \emph{From now on, all
 the graphs we  consider are connected}.\\

Let us now recall the relations of induction satisfied by the Tutte
polynomial \cite{Tutte:dichromate}. (These relations reminiscent of the relations of induction of the chromatic polynomial \cite{Whitney:chromatic} are easy to prove from \Ref{eq:Tutte-explicit}).
%
%
\begin{prop}\textbf{(Tutte)} \label{prop:induction}
Let $G$ be a graph and $e$ be any edge of $G$. The Tutte polynomial of $G$ satisfies:
\begin{eqnarray} \label{eq:induction}
\begin{array}{l|ll}
T_{G}(x,y)~=
&~x\cdot T_{G{/e}}(x,y) & \textrm{ if $e$ is an isthmus,} \\
&~y\cdot T_{G{\backslash e}}(x,y) & \textrm{ if $e$ is a loop,} \\
&~T_{G{/e}}(x,y)~+~T_{G{\backslash e}}(x,y) & \textrm{ if $e$ is neither a loop nor an isthmus.}\\
\end{array}
\end{eqnarray}
\end{prop}

\vspace{.5cm}

Before we close this section, we present another characterization of
the Tutte polynomial that will have a great similarity with ours.
This characterization, given by Tutte in 1954, uses
 the notion of \emph{activity}.  Let $G=(V,E)$ be
a graph. Suppose that the edge set $E$ is linearly ordered. Given a
spanning tree $T$, an external (resp. internal) edge is said
\emph{active} if it is minimal in its fundamental cycle
(resp. cocycle). We recall the following theorem by Tutte
\cite{Tutte:dichromate}. 

\begin{thm} \textbf{(Tutte)} \label{thm:activity-Tutte}
Let $G=(V,E)$ be a connected graph. Given any linear order on the edge set $E$, the Tutte polynomial of $G$ is equal to 
\begin{eqnarray} \label{eq:sum-activity-Tutte}
T_G(x,y)=\sum_{T \textrm{ spanning tree}} x^{i(T)}y^{e(T)},
\end{eqnarray}
where the sum is over all spanning trees and $i(T)$ (resp. $e(T)$) is the number of active internal (resp. external) edges.
\end{thm}

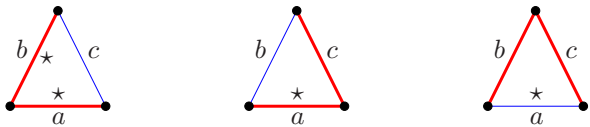
\begin{figure}[htb!]
\begin{center}
\input{exp-Tutte-triangle2.pstex_t}
\caption{The spanning trees of $K_3$.}\label{fig:exp-Tutte-triangle}
\end{center}
\end{figure}

For instance, if $G$ is $K_3$ there are three spanning trees that we represented in Figure \ref{fig:exp-Tutte-triangle}. If the linear order  on the edge set is $a<b<c$ then the active edges are the one indicated by a $\star$. Hence, the spanning trees (taken from left to right) have respective contributions $x^2$, $x$ and $y$ and the Tutte polynomial is $T_{K_3}(x,y)=x^2+x+y$.\\

Theorem \ref{thm:activity-Tutte} is  surprising because it implies
 that the sum in \Ref{eq:sum-activity-Tutte} does not depend on the
 ordering of the edge set (whereas the summands clearly depends
 on that order). However,  this theorem is easily proved
 by induction. Indeed, it is simple to prove that
 the induction relation of Proposition \ref{prop:induction} holds 
 \emph{for the  edge having the largest label}.\\



\section{The Tutte polynomial of embedded graphs}\label{section:def-Tutte-OB}
In this section, we present a new characterization of the Tutte
polynomial. This  characterization is, just as in Theorem
\ref{thm:activity-Tutte}, a generating function of the spanning trees counted according to activities 
(see Theorem \ref{thm:activity-OB}). We will only change the
definition of activities, using the embedding structure 
%
instead of a linear order on the edge set.\\

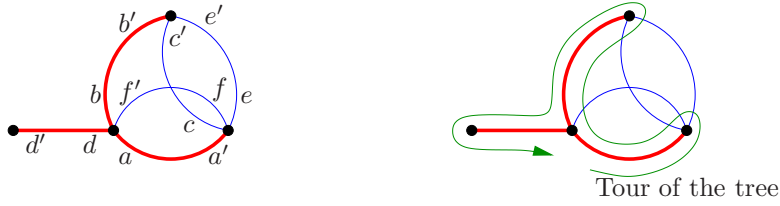
\begin{figure}[htb!]
\begin{center}
\input{tour-of-tree.pstex_t}
\caption{Intuitive representation of the tour of a spanning tree (indicated by thick lines).}\label{fig:tour-of-tree}
\end{center}
\end{figure}

Let us begin by defining the \emph{tour of a spanning
  tree}. Informally, the \emph{tour of a tree} is 
a walk around the tree 
that  follows internal edges and crosses external edges. 
A graphical representation of this notion is represented in
Figure \ref{fig:tour-of-tree}. \\ 

Let $\mathcal{G}=(H,\sigma,\alpha)$ be an embedding of the graph $G=(V,E)$. Given a spanning tree $T$, we define a function $t$ on the set $H$ of half-edges by :
\begin{eqnarray} \label{def:t}
\begin{array}{l|ll}
t(h)=& \sigma(h) & \textrm{ if } h  \textrm{ is external,}\\
& \sigma\ccirc \alpha(h) & \textrm{ if } h \textrm{ is internal.}
\end{array}
\end{eqnarray}
The mapping $t$ is called the \emph{motion function}. The motion
function $t$ is a bijection on $H$ (since the inverse mapping is easily seen
to be $t^{-1}(h)=\sigma^{-1}(h)$ if $\sigma^{-1}(h)$ is external and
$t^{-1}(h)=\alpha\ccirc \sigma^{-1}(h)$ if  $\sigma^{-1}(h)$ is
internal). In fact, as we will prove shortly (Lemma
\ref{thm:tour-is-cyclic}), the motion function $t$ is  a cyclic
permutation. For instance, for the embedded graph of Figure
\ref{fig:tour-of-tree}, the motion function is the cycle
$(a,e,f,c,a',f',b,c',e',b',d,d')$. The cyclic order defined by the
motion function $t$ on the set of half-edges is what we call the
\emph{tour of the tree} $T$. Note that turning around a spanning tree
and writing down the half-edges in order of appearance gives an
encoding of the tree. This encoding is closely related to the
encodings of maps with a distinguished spanning tree given 
%
by Lehman and Walsh \cite{Walsh:counting-maps-2}. \\

Before we prove that the motion function $t$ is a cyclic permutation, we describe what happens to this function when an edge is deleted or contracted.

\begin{lemma}\label{thm:induction-tour}
Let $\mathcal{G}=(H,\sigma,\alpha)$ be an embedded graph, $T$ 
 a spanning tree and $t$  the corresponding motion function. For all external (resp. internal)
 edge $e=\{h_1,h_2\}$, the spanning tree $T$ (resp. $T-{e}$) of
 $\mathcal{G}_{\backslash e}$ (resp. $\mathcal{G}_{/e}$) defines a
 motion function $t'$ on $H\setminus\{h_1,h_2\}$ such that  
%
%
$$\begin{array}{l|ll}
t'(h)=
& t\ccirc t\ccirc t(h) &\textrm{if } (t(h)=h_1 \textrm{ and } t(h_1)=h_2) \textrm{ or } (t(h)=h_2 \textrm{ and } t(h_2)=h_1), \\
& t\ccirc t(h) & \textrm{if } (t(h)=h_1 \textrm{ and } t(h_1)\neq h_2) \textrm{ or } (t(h)=h_2 \textrm{ and } t(h_2)\neq h_1), \\
&t(h) & \textrm{otherwise}. 
\end{array}
$$
\end{lemma}

\dem Lemma \ref{thm:induction-tour} follows immediately from the definitions and Equations  \Ref{eq:deletion-plongement} and \Ref{eq:contraction-plongement}.
\findem

\noindent \textbf{Remark:} Another way of stating Lemma \ref{thm:induction-tour} is to say that the cycles of the permutation $t'$ are obtained from the cycles of $t$ by erasing  $h_1$ and $h_2$. 
Consider, for instance, the embedded graph  and the spanning tree represented in Figure \ref{fig:tour-of-tree}. The motion function is the cycle $t=(a,e,f,c,a',f',b,c',e',b',d,d')$. If we delete the edge the external edge  $\{e,e'\}$ (resp. internal edge  $\{b,b'\}$), the motion function becomes $t'=(a,f,c,a',f',b,c',b',d,d')$ (resp. $t'=(a,e,f,c,a',f',c',e',d,d')$).\\

\begin{lemma}\label{thm:tour-is-cyclic}
For any embedded graph and any spanning tree, the motion function is a cyclic permutation.
\end{lemma}

\dem We prove the lemma by induction on the number of edges of the graph. The property is obviously true for the graph reduced to a loop and the graph reduced to an isthmus.
We assume the property holds for all graphs with at most $n\geq 1$
edges and consider an embedded graph $\mathcal{G}$ with $n+1$
edges. Let $T$ be a spanning tree and $t$ the corresponding motion function. We know that $t$ is a permutation, that is, a product of cycles. We consider an edge $e=\{h_1,h_2\}$. \\
\ite \emph{In any cycle of the motion function $t$, there is a half-edge $h \neq h_1, h_2$.}\\
First note that $t(h_i)\neq h_i$
for $i=1,2$. Indeed, if $e$ is external, this would mean
$\sigma(h_i)=h_i$ which is excluded or $e$ would be an isthmus not in
the spanning tree. Similarly, if  $e$ is internal, we would have
$\sigma\ccirc \alpha(h_i)=h_i$ which is excluded or $e$ would be a loop
in the spanning tree. Moreover, we cannot have ($t(h_1)=h_2$ and
$t(h_2)=h_1$). Indeed, this would mean that $e$ is either an isolated
loop (if $e$ is external) or an isolated isthmus (if $e$ is internal)
contradicting our hypothesis that $\mG$ is connected and has more than one edge. \\ 
\ite \emph{The motion function is cyclic.}\\
If $e$ is external (resp. internal), we consider the motion function
$t'$ defined by the spanning tree $T$ (resp. $T-e$) on
$\mG_{\backslash e}$ (resp. $\mG_{/e}$). By Lemma
\ref{thm:induction-tour},  the cycles of $t'$ are the cycles of $t$ where the half-edges
 $h_1,h_2$ are erased. 
%
Suppose now that  $t$ is not cyclic.  Then $t$ has at least two cycles each containing a half-edge $h \neq h_1, h_2$. Therefore, $t'$ has at least two non-empty cycles, which contradicts our induction hypothesis.
\findem

Consider a rooted embedded graph $\mathcal{G}=(H,\sigma,\alpha,h)$ and
a spanning tree $T$. We define a linear order on $H$ by
$h<t(h)<t^2(h)<\ldots<t^{|H|-1}(h)$, where $t$ is the motion
function. (This defines a linear order on the set of half-edges since
$t$ is a cyclic permutation.) We call \emph{$(\mG,T)$-order} this
order. We define a  $(\mG,T)$-order on 
 the edge set by setting
$e=\{h_1,h_2\}<e'=\{h_1',h_2'\}$ if $\min(h_1,h_2)<\min(h_1',h_2')$. Note that this is also a linear order.\\ 

We are now ready to define a new notion of activity (that we call \emph{embedding-activity}, so as to differentiate it from Tutte's notion of activity) and state our main theorem.

\begin{Def}
Let  $\mG$ be a rooted embedded graph and $T$ be a spanning tree.   
We say that an external (resp. internal) edge is
\emph{$(\mG,T)$-active} (or \emph{embedding-active} if $\mG$ and $T$
are clear from the context) if it is minimal for the $(\mG,T)$-order in its fundamental cycle (resp. cocycle).  
\end{Def}

\noindent \textbf{Example:} Consider the embedded graph $\mG$ and spanning tree $T$ represented in Figure \ref{fig:tour-of-tree}. Take the half-edge $a$ as the root of  $\mG$. The $(\mG,T)$-order on the half-edges is $a<e<f<c<a'<f'<b<c'<e'<b'<d<d'$. Therefore, the $(\mG,T)$-order on the edges is $\{a,a'\}<\{e,e'\}<\{f,f'\}<\{c,c'\}<\{b,b'\}<\{d,d'\}$. The internal active edges are $\{a,a'\}$ and $\{d,d'\}$ and there is no external active edge. For instance, $\{e,e'\}$ is not active since $\{a,a'\}$ is in its fundamental cycle. \\

\begin{thm} \label{thm:activity-OB}
Let $\mG$ be any rooted embedding of the connected graph $G$ (with at least one edge).  The Tutte polynomial of $G$ is equal to 
\begin{eqnarray}
T_G(x,y)=\sum_{T \textrm{ spanning tree}} x^{\mI(T)}y^{\mE(T)}, \label{eq:Tutte-embedded}
\end{eqnarray}
where the sum is over all spanning trees and $\mI(T)$ (resp. $\mE(T)$) is the number of embedding-active internal (resp. external) edges.
\end{thm}

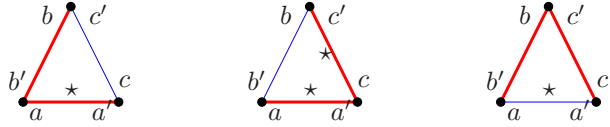
\begin{figure}[htb!]
\begin{center}
\input{exp-Tutte-triangle3.pstex_t}
\caption{The embedding-activities of the spanning trees of $K_3$.}\label{fig:exp-Tutte-triangle-OB}
\end{center}
\end{figure}

We represented the spanning trees of $K_3$ in Figure \ref{fig:exp-Tutte-triangle-OB}. If the embedding is rooted on the half-edge $a$, then the embedding-active edges are the one indicated by a $\star$. Hence, the spanning trees (taken from left to right) have respective contributions $x$, $x^2$ and $y$ and the Tutte polynomial is $T_{K_3}(x,y)=x^2+x+y$.\\

Let us  emphasize that Theorem \ref{thm:activity-OB} is \emph{not}
a special case of Theorem \ref{thm:activity-Tutte}. Indeed, the
$(\mG,T)$-order  defined above is a linear order on the edge set that \emph{depends}
 on the tree $T$.\\

The rest of this section is devoted to the proof of Theorem \ref{thm:activity-OB}. 

\begin{lemma} \label{thm:induction-T-order}
Let $\mG$ be a rooted embedded graph with edge set $E$ and half-edge set $H$. Let $T$ be a spanning tree and $e=\{h_1,h_2\}$ be an edge not containing the root. If $e$ is external (resp. internal), the $(\mG_{\backslash e},T)$-order (resp. $(\mG_{/ e},T-e)$-order) on  $H\setminus\{h_1,h_2\}$ and $E-e$ is simply the restriction of the  $(\mG,T)$-order to these sets.\\
\end{lemma}

\noindent \textbf{Proof of Lemma \ref{thm:induction-T-order}: }
By Lemma \ref{thm:induction-tour}, we see that if $e$ is external
(resp. internal), the $(\mG_{\backslash e},T)$-order (resp. $(\mG_{/
  e},T-e)$-order) on  the half-edge set $H\setminus\{h_1,h_2\}$ is simply the
restriction of the  $(\mG,T)$-order to this set. The same property follows immediately for the edge set.
\findem

\noindent \textbf{Proof of Theorem \ref{thm:activity-OB}: }
We associate to the rooted embedded graph $\mathcal{G}$ the polynomial 
$$\mT_{\mG}(x,y)=\sum_{T \textrm{ spanning tree}} x^{\mI(T)}y^{\mE(T)},$$
where  $\mI(T)$ (resp. $\mE(T)$) is the number of embedding-active
internal (resp. external) edges. We want to show that 
the polynomial $\mT_{\mG}(x,y)$ is equal to the Tutte polynomial
$T_{G}(x,y)$ of 
$G$. We proceed by induction on the number of edges, using Proposition \ref{prop:induction}.\\  
\ite 
The graphs with one edge are the graph $L$ reduced to a loop and the graph $I$ reduced to an isthmus. The graph $L$ (resp. $I$) has a unique rooted embedding $\mathcal{L}$ (resp. $\mathcal{I}$).  We check that 
$T_L(x,y)=y=\mathcal{T}_\mathcal{L}(x,y)$ and $T_I(x,y)=x=\mathcal{T}_\mI(x,y)$.\\
\ite We assume the property holds for all (connected) graphs with at most $n\geq 1$ edges and consider a rooted embedding $\mathcal{G}=(H,\sigma,\alpha,h_0)$ of a graph $G$ with $n+1$ edges. We denote by $v_0$ the vertex incident to the root $h_0$  and $e_0$ the edge containing $h_0$. We denote by $h_*=\sigma^{-1}(h_0)$ the half-edge preceding $h_0$ around $v_0$ and by $e_*=\{h_*,h_*'\}$ the edge containing $h_*$. \\
We study separately the 3 different cases of the induction relation \Ref{eq:induction}.\\
\noindent \emph{\textbf{Case 1:} The edge $e_*$ is neither an isthmus nor a loop.}\\
The set $\mathbb{T}$ of spanning trees of $G$ can be partitioned into
$\mathbb{T}=\mathbb{T}_1\cup \mathbb{T}_2$, where $\mathbb{T}_1$
(resp. $\mathbb{T}_2$) is the set of spanning trees containing
(resp. not containing) the edge $e_*$. The set $\mathbb{T}_1$
(resp. $\mathbb{T}_2$) is in bijection by the mapping
$\Phi_1:~T\mapsto T-e_*$ (resp. $\Phi_2:~T\mapsto T$) with the
spanning trees of $G_{/e_*}$ (resp. $G_{\backslash e_*}$). We want to
show $e_*$ is never embedding-active and that the mappings $\Phi_i$ preserve the
embedding-activities: for any tree $T$ in $\mathbb{T}_1$
(resp. $\mathbb{T}_2$), an edge is $(\mG,T)$-active if and only if it
is $(\mG_{/ e_*},T-e_*)$-active (resp. $(\mG_{\backslash
  e_*},T)$-active). We are going to prove successively the following four points:
\begin{itemize}
\item \emph{The edges $e_*$ and $e_0$ are distinct.} \\ 
First note that $h_0\neq h_*$ or we would have  $\sigma(h_*)=h_*$ implying that  $v_0$ has degree one hence that $e_*$ is an isthmus. 
Also,  $h_0\neq h_*'$ or we would have $\sigma(h_*)=\alpha(h_*)$ implying that $e_*$ is a loop. Thus, $e_*=\{h_*,h_*'\}$ does not contain $h_0$.
\item \emph{Given any spanning tree, the edge $e_*$ is  maximal in its fundamental cycle or cocycle.} \\
Let $T$ be a spanning tree of $G$. 
Suppose first that the edge $e_*$ is internal. In this case, the motion function $t$ satisfies,
$t(h_*')=\sigma\alpha(h_*')=h_0$. Hence, $h_*'$ is the greatest half-edge for
the  $(\mG,T)$-order.  Let $C_*$ be the fundamental cocycle of
$e_*$. Recall that, apart from $e_*$, the edges in $C_*$ are
external.  Let $e$ be one of those edges. We want to prove that $e<e_*$ for the $(\mG,T)$-order. 
Consider the embedded graph $\mG'$ (with two edges) obtained by deleting all external edges distinct from $e$ and  contracting all internal edges distinct from $e_*$.  
We apply the same operations to the tree $T$, thus obtaining $T'=\{e_*\}$.
By induction, we see that in $(\mG',T')$ the edge $e$ is in the
fundamental cocycle of $e_*$. Hence, the only possibility for $\mG'$ is
the embedded graph represented  on the left of Figure
\ref{fig:case-study-external}. The motion function of  $(\mG',T')$ is
 the cycle $(h,h_*,h',h_*')$ where $h,h'$ are the half-edges of
$e$. By (repeated application of) Lemma \ref{thm:induction-tour}, the
half-edges $h,h_*,h',h_*'$  appear in the same cyclic order around the
spanning tree $T$ of  $\mG$. Moreover, since $h_*'$ is the greatest
half-edge for the  $(\mG,T)$-order, we have $h<h_*<h'<h_*'$ hence $e<e_*$ for the $(\mG,T)$-order. 
Thus, $e_*$ is maximal (for the $(\mG,T)$-order) in its fundamental cocycle. \\ 
Similarly, if $e_*$ is external,  $h_*$ is the greatest half-edge for the  $(\mG,T)$-order (since
$t(h_*)=\sigma(h_*)=h_0$). Then, we consider the fundamental cycle $C_*$ of
$e_*$. For any edge $e$ in $C_*-e_*$, we consider the embedded graph $\mG'$  obtained by deleting all
external edges distinct from $e_*$ and  contracting all internal edges
distinct from $e$. Applying the same operations to the tree $T$, we
obtain $T'=\{e\}$. We see that the only possibility for $\mG'$ is
the embedded graph represented  on the right of Figure
\ref{fig:case-study-external}. The motion function of  $(\mG',T')$ is
the cycle $(h,h_*',h',h_*)$ where $h,h'$ are the half-edges of
$e$. We conclude as above.\\ 
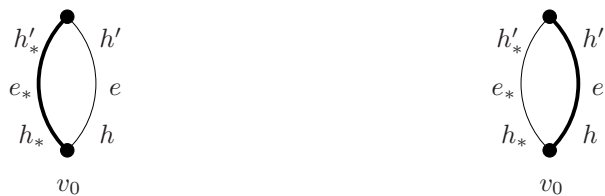
\begin{figure}[htb!]
\begin{center}
\input{case-study-external2.pstex_t}
\caption{The only possibility for $\mG'$ when $e_*$ is internal (left)
  or external (right).}\label{fig:case-study-external}
\end{center}
\end{figure}
\item \emph{For any tree $T$ in $\mathbb{T}_1$ (resp. $\mathbb{T}_2$), the  $(\mG,T)$-active and  $(\mG_{/ e_*},T-e_*)$-active  (resp. $(\mG_{\backslash e_*},T)$-active) edges are the same.}\\ 
First note that $e_*$ is never alone in its fundamental cycle or cocycle (or $e_*$ would be a loop or isthmus). Hence, by the preceding point, $e_*$ is never embedding-active. We now look at the embedding-activities of the other edges. Let  $T$ be a tree in $\mathbb{T}_1$ (i.e. containing $e_*$). 
Let $e$ be an external (resp. internal) edge distinct from $e_*$ and 
let $C$ be its fundamental cycle (resp. cocycle).  The fundamental cycle
(resp. cocycle) of $e$ in $(\mG_{/ e_*},T-e_*)$ is $C-e_*$.  Note that
the $(\mG,T)$-minimal element of $C$ is in $C-e_*$ (since, if $e_*$ is
in $C$ then $e$ is in the fundamental cycle of $e_*$ hence $e<e_*$ for
the $(\mG,T)$-order). Moreover, by Lemma \ref{thm:induction-T-order},
the  $(\mG,T)$-order and $(\mG_{/ e_*},T-e_*)$-order coincide on
$C-e_*$.  Hence, the  $(\mG,T)$-minimal element of $C$ is the $(\mG_{/
  e_*},T-e_*)$-minimal element in $C-e_*$. Therefore, the edge $e$ is
$(\mG,T)$-active if and only if it is $(\mG_{/ e_*},T-e_*)$-active.\\ 
The case where $T$ is a tree in $\mathbb{T}_2$ (i.e. not containing $e_*$) is identical.
\item \emph{The polynomial $\mT_\mG(x,y)$ is equal to the Tutte polynomial  $T_G(x,y)$.}\\
From the properties above, we have 
\begin{eqnarray}
\mT_\mG(x,y)&\displaystyle \equiv& \sum_{T \textrm{ spanning tree of } G} x^{\mI(T)}y^{\mE(T)}\nonumber\\
& \displaystyle =&\sum_{T \in \mathbb{T}_1} x^{\mI(T)}y^{\mE(T)}+ \sum_{T \in \mathbb{T}_2} x^{\mI(T)}y^{\mE(T)}\nonumber\\
& \displaystyle =&\sum_{T \in \mathbb{T}_1} x^{\mI'(T-e_*)}y^{\mE'(T-e_*)}+\sum_{T \in \mathbb{T}_2} x^{\mI''(T)}y^{\mE''(T)}\label{eq:partition-trees}
\end{eqnarray}
where $\mI'(T-e_*)$, $\mE'(T-e_*)$, $\mI''(T)$, $\mI''(T)$ are
respectively the number of internal $(\mG_{/e_*},T-e_*)$-active,
external $(\mG_{/e_*},T-e_*)$-active,  $(\mG_{\backslash
  e_*},T)$-active  and external $(\mG_{\backslash e_*},T)$-active
edges. \\
%
%
In the right-hand side of \Ref{eq:partition-trees} we recognize the polynomials
 $\mT_{\mG_{/e_*}}(x,y)$ and $\mT_{\mG_{\backslash e_*}}(x,y)$. By 
the induction hypothesis,  these polynomials are the Tutte polynomials  $T_{G_{/e_*}}(x,y)$ and $T_{G_{\backslash e_*}}(x,y)$. Thus, 
\begin{eqnarray}\label{eq:partition-trees2}
\hspace{1cm} \mT_\mG(x,y)~=~\mT_{\mG_{/e_*}}(x,y)+\mT_{\mG_{\backslash e_*}}(x,y)=  T_{G_{/e_*}}(x,y)+T_{G_{\backslash e_*}}(x,y).
\end{eqnarray}
In view of the induction relation of Proposition \ref{prop:induction}, this is the Tutte polynomial $T_G(x,y)$. 
\end{itemize} 
\vspace{.5cm}

\noindent \emph{\textbf{Case 2}: The edge $e_*$ is an isthmus.}\\
Since $e_*$ is an isthmus, it is in every spanning tree. Moreover, being alone in its fundamental cocycle, it is always active. We want to show that for any spanning tree $T$, the embedding-activity of any edge other than $e_*$ is the same in  $(\mG,T)$ and in  $(\mG_{/e_*},T-e_*)$. Before we do that, we must cope with a (little) technical difficulty: the edge $e_*$ might be equal to $e_0$ in which case we should specify how to root the graph  $\mG_{/e_*}$. \\
First note that  $h_0\neq h_*'$ or we would have $\sigma(h_*)=\alpha(h_*)$ implying that $e_*$ is a loop. Suppose now that  $h_0= h_*$ (equivalently, $\sigma(h_*)=h_*$). In this case, we define the root of $\mG_{/e_*}$ to be $h_1=\sigma(h_*')$ ($h_1$ is not an half-edge of $e_*$ or $e_*$ would be an isolated isthmus).
\begin{itemize}
\item \emph{For any spanning tree $T$ of $G$, the  $(\mG,T)$-order and the  $(\mG_{/e_*},T-e_*)$-order coincide on $E-e_*$.}\\
If $e_*\neq e_0$ the property is given by Lemma \ref{thm:induction-T-order}. Now suppose that $e_*= e_0$ (that is $h_*=h_0$). Since $e_*$ is internal, the motion function $t$ satisfies $t(h_*')=h_0$ and $t(h_0)=h_1$. Therefore, the $(\mG,T)$-order on half-edges is $h_0=h_*<h_1<t(h_1)<\ldots<h_*'$. Let us denote by $\mG_1$ the embedded graph $\mG$ rooted on $h_1$.  The $(\mG_1,T)$-order on half-edges is $h_1<t(h_1)<\ldots<h_*'<h_0=h_*$. Thus, the  $(\mG_1,T)$-order and $(\mG,T)$-order coincide on $E-e_*$. Moreover, by Lemma \ref{thm:induction-T-order}, the $(\mG_1,T)$-order and $(\mG_{/e_*},T)$-order   coincide on  $E-e_*$.
\item \emph{For any spanning tree $T$, the set of $(\mG,T)$-active edges distinct from $e_*$ is the set of  $(\mG_{/ e_*},T-e_*)$-active edges.}\\
For any tree $T$ and any external (resp. internal) edge $e\neq e_*$, the fundamental cycle (resp. cocycle) of $e$  does not contain $e_*$ and is the same in  $(\mG,T)$ and in  $(\mG_{/e_*},T-e_*)$. Since the   $(\mG,T)$-order and the $(\mG_{/e_*},T-e_*)$-order coincide on $E-e_*$, the edge $e$ is  $(\mG,T)$-active if and only if it is  $(\mG_{/e_*},T-e_*)$-active.
\item  \emph{The polynomial $\mT_\mG(x,y)$ is equal to the Tutte polynomial  $T_G(x,y)$.}\\
From the properties above, we have 
\begin{eqnarray}
\mT_\mG(x,y)&  \displaystyle  \equiv& \sum_{T \textrm{ spanning tree of } G} x^{\mI(T)}y^{\mE(T)} \nonumber\\
& \displaystyle =& \sum_{T \textrm{ spanning tree of } G} x^{1+\mI'(T-e_*)} y^{\mE'(T-e_*)}\nonumber\\
& \displaystyle =& x \cdot  \sum_{T \textrm{ spanning tree of } G} x^{\mI'(T-e_*)}y^{\mE'(T-e_*)}\label{eq:partition-trees-isthmus}
\end{eqnarray}
where $\mI'(T-e_*)$ and  $\mE'(T-e_*)$ are respectively the number of internal $(\mG_{/e_*},T-e_*)$-active and external $(\mG_{/e_*},T-e_*)$-active edges. \\
In the right-hand side of \Ref{eq:partition-trees-isthmus} we recognize the sum as being $\mT_{\mG_{/e_*}}(x,y)$.
By the induction hypothesis, we know this polynomial to be equal to the Tutte polynomial  $T_{G_{/e_*}}(x,y)$. Thus,
\begin{eqnarray}\label{eq:partition-trees2-isthmus}
\mT_\mG(x,y)~=~x\cdot \mT_{\mG_{/e_*}}(x,y)=  x\cdot T_{G_{/e_*}}(x,y).
\end{eqnarray}
In view of the induction relation of Proposition \ref{prop:induction}, this is the Tutte polynomial $T_G(x,y)$. 
\end{itemize}
\vspace{.5cm}

\noindent \emph{\textbf{Case 3:} The edge $e_*$ is a loop.}\\
This case is dual to Case 2. \\
Since $e_*$ is a loop, it is always external and always active.  We
want to show that for any spanning tree $T$, the embedding-activity  
of any edge other than $e_*$
is the same in  $(\mG,T)$ and in  $(\mG_{\backslash e_*},T)$. Before
we do that, we must choose a root for  $\mG_{\backslash e_*}$ when
$e_*=e_0$. We see that $h_0\neq h_*$ or we would have
$\sigma(h_*)=h_*$ implying that $e_*$ is an isthmus. Suppose now that
$h_0= h_*'$ (equivalently, $\alpha(h_*)=\sigma(h_*)$). In this
case, we define the root of  $\mG_{\backslash e_*}$ to be $h_1=\sigma(h_0)$ ($h_1$ is not an half-edge of $e_*$ or $e_*$ would be an isolated loop). 
\begin{itemize}
\item \emph{For any spanning tree $T$ of $G$, the  $(\mG,T)$-order  and the $(\mG_{\backslash e_*},T)$-order coincide on $E-e_*$.}\\
The proof of Case 2 can be copied verbatim except  ``$e_*$ is internal'' is replaced by ``$e_*$ is external''.
\item \emph{For any spanning tree $T$, the set of $(\mG,T)$-active edges distinct from $e_*$ is the set of  $(\mG_{\backslash e_*},T)$-active edges.}\\
The proof of Case 2 can be copied verbatim.
\item  \emph{The polynomial $\mT_\mG(x,y)$ is equal to the Tutte polynomial  $T_G(x,y)$.}\\
From the properties above, we have 
\begin{eqnarray}
\mT_\mG(x,y)& \displaystyle \equiv& \sum_{T \textrm{ spanning tree of } G} x^{\mI(T)}y^{\mE(T)}\nonumber\\
&\displaystyle  =& y\cdot \sum_{T \textrm{ spanning tree of } G} x^{\mI''(T)}y^{\mE''(T)}\label{eq:partition-trees-loop}
\end{eqnarray}
where $\mI''(T)$ and  $\mE''(T)$ are respectively the number of internal $(\mG_{\backslash e_*},T)$-active and external $(\mG_{\backslash e_*},T)$-active edges. \\
In the right-hand side of \Ref{eq:partition-trees-loop} we recognize the sum as being $\mT_{\mG_{\backslash e_*}}(x,y)$. 
By the induction hypothesis, we know this polynomial to be equal to the Tutte polynomial  $T_{G_{\backslash e_*}}(x,y)$. Thus,
\begin{eqnarray}\label{eq:partition-trees2-loop}
\mT_\mG(x,y)~=~y\cdot \mT_{\mG_{\backslash e_*}}(x,y)=  y\cdot T_{G_{\backslash e_*}}(x,y).
\end{eqnarray}
In view of the induction relation of Proposition \ref{prop:induction}, this is the Tutte polynomial $T_G(x,y)$. 
\end{itemize}
~ \findem

\section{Concluding remarks}\label{section:remarks}
We conclude this paper by mentioning some possible
applications of our characterization of the Tutte polynomial. 
Some of them will be developed in a forthcoming
paper~\cite{OB:Tutte-plongement}. \\

The characterization of the Tutte polynomial in terms   of the
 activities of spanning trees is sometimes thought of as slightly 
unnatural. 
It is true that the dependence of this characterization on
a particular linear ordering of the edge set is a bit puzzling. We want to argue that an embedding may be a less arbitrary 
 structure than a linear order on the edge
set. As a matter of fact, there are a number of problems in which the
embedding structure is explicitly given. \\ 
Firstly, some famous conjectures 
deal with the Tutte polynomial, or sometimes the chromatic polynomial, of
 \emph{planar graphs}. A graph is \emph{planar} if and only if can be
 drawn on the plane without intersecting edges. Equivalently,  it
 has an embedding $(H,\sigma,\alpha)$ with Euler characteristic equal
 to 2.
For instance, the four color theorem can be stated as: \emph{$\mT_\mG(-3,0)\neq 0$ for any loopless planar embedding $\mG$}.\\ 
Another type of problem where the embedding structure is explicitly
 given appears in mathematical physics, in the study of the \emph{Potts model on random lattices}
 \cite{Baxter:Potts-model,Bonnet:Potts-loop-equations,Daul:Potts}. The
 \emph{Potts model} is an important statistical mechanics model for
 particles interacting in a discrete space (i.e. a graph)
 \cite{Baxter:exactly-solved-model,Sokal-multivariate-Tutte}. It was
 shown by Fortuin and Kasteleyn \cite{Fortuin:Tutte=Potts} that the
 partition function of the Potts model on the graph $G$ is equivalent
 (up to a change of variables) to the Tutte polynomial of $G$. Studying
 this model on a \emph{random lattice} means that the underlying space
 (the graph) is random. Usually, the underlying space is
 supposed to have a uniform distribution over a class of rooted
 maps.  In this case, the partition
 function of the model is equivalent to the sum of the Tutte
 polynomial over the class of maps. For instance, if the \emph{random
 lattice} is understood as the uniform distribution on the set of
 rooted planar maps with $n$ edges, we are to study the partition
 function 
$$Z_n(x,y)=\sum_{\mG \textrm{ rooted map with $n$ edges}} \mT_\mG(x,y).$$
Note that this partition function can also be written as
$$Z_n(x,y)=\sum_{\mG,~T} x^{\mI(\mG,T)}y^{\mE(\mG,T)},$$
where the sum is over all rooted maps $\mG$ with $n$ edges and all spanning trees $T$  and $\mI(\mG,T)$ (resp. $\mE(\mG,T)$) is the number of $(\mG,T)$-active internal (resp. external) edges. This last expression of the partition function could be interesting since there are several nice encodings for \emph{rooted maps with a distinguished spanning tree} (a.k.a. \emph{tree-rooted maps}) \cite{OB:Boisees,Walsh:counting-maps-2,Mullin:tree-rooted-maps}.\\

Let us also mention that embedding structures can be used to define
several bijections between spanning trees and some other structures counted by the Tutte polynomial, including root-connected out-degree sequences and recurrent configurations of the sandpile model \cite{OB:Tutte-plongement}. Our characterization of the Tutte polynomial is deeply related to these bijections and gives a combinatorial interpretation of several enumerative results.\\

\textbf{Acknowledgments: } This work has benefited from discussions
with Yvan Le Borgne, Emeric Gioan and Michel Las Vergnas. I would also
like to thank my advisor Mireille Bousquet-Mélou for her constant
support and guidance.

\bibliography{../../../biblio/allref}
\bibliographystyle{plain}

\end{document}

%% file: exp-embedding.pstex_t
\begin{picture}(0,0)%
\epsfig{file=exp-embedding.pstex}%
\end{picture}%
\setlength{\unitlength}{1381sp}%
\begingroup\makeatletter\ifx\SetFigFont\undefined%
\gdef\SetFigFont#1#2#3#4#5{%
  \reset@font\fontsize{#1}{#2pt}%
  \fontfamily{#3}\fontseries{#4}\fontshape{#5}%
  \selectfont}%
\fi\endgroup%
\begin{picture}(15419,4020)(-2518,-5236)
\put(3001,-2761){\makebox(0,0)[lb]{\smash{\SetFigFont{5}{6.0}{\rmdefault}{\mddefault}{\updefault}\normalsize{$c'$}}}}
\put(3301,-4561){\makebox(0,0)[lb]{\smash{\SetFigFont{5}{6.0}{\rmdefault}{\mddefault}{\updefault}\normalsize{$c$}}}}
\put(3826,-5236){\makebox(0,0)[lb]{\smash{\SetFigFont{5}{6.0}{\rmdefault}{\mddefault}{\updefault}\normalsize{$a'$}}}}
\put(2026,-5236){\makebox(0,0)[lb]{\smash{\SetFigFont{5}{6.0}{\rmdefault}{\mddefault}{\updefault}\normalsize{$a$}}}}
\put(1951,-2461){\makebox(0,0)[lb]{\smash{\SetFigFont{5}{6.0}{\rmdefault}{\mddefault}{\updefault}\normalsize{$b'$}}}}
\put(3751,-2311){\makebox(0,0)[lb]{\smash{\SetFigFont{5}{6.0}{\rmdefault}{\mddefault}{\updefault}\normalsize{$e'$}}}}
\put(1201,-5011){\makebox(0,0)[lb]{\smash{\SetFigFont{5}{6.0}{\rmdefault}{\mddefault}{\updefault}\normalsize{$d$}}}}
\put(  1,-5011){\makebox(0,0)[lb]{\smash{\SetFigFont{5}{6.0}{\rmdefault}{\mddefault}{\updefault}\normalsize{$d'$}}}}
\put(4501,-3961){\makebox(0,0)[lb]{\smash{\SetFigFont{5}{6.0}{\rmdefault}{\mddefault}{\updefault}\normalsize{$e$}}}}
\put(1351,-3961){\makebox(0,0)[lb]{\smash{\SetFigFont{5}{6.0}{\rmdefault}{\mddefault}{\updefault}\normalsize{$b$}}}}
\put(1951,-3886){\makebox(0,0)[lb]{\smash{\SetFigFont{5}{6.0}{\rmdefault}{\mddefault}{\updefault}\normalsize{$f'$}}}}
\put(3901,-3886){\makebox(0,0)[lb]{\smash{\SetFigFont{5}{6.0}{\rmdefault}{\mddefault}{\updefault}\normalsize{$f$}}}}
\put(8401,-5011){\makebox(0,0)[lb]{\smash{\SetFigFont{5}{6.0}{\rmdefault}{\mddefault}{\updefault}\normalsize{$d'$}}}}
\put(12901,-3961){\makebox(0,0)[lb]{\smash{\SetFigFont{5}{6.0}{\rmdefault}{\mddefault}{\updefault}\normalsize{$e$}}}}
\put(9751,-3961){\makebox(0,0)[lb]{\smash{\SetFigFont{5}{6.0}{\rmdefault}{\mddefault}{\updefault}\normalsize{$b$}}}}
\put(11701,-2911){\makebox(0,0)[lb]{\smash{\SetFigFont{5}{6.0}{\rmdefault}{\mddefault}{\updefault}\normalsize{$c'$}}}}
\put(10501,-4261){\makebox(0,0)[lb]{\smash{\SetFigFont{5}{6.0}{\rmdefault}{\mddefault}{\updefault}\normalsize{$f'$}}}}
\put(11701,-4261){\makebox(0,0)[lb]{\smash{\SetFigFont{5}{6.0}{\rmdefault}{\mddefault}{\updefault}\normalsize{$f$}}}}
\put(12226,-3886){\makebox(0,0)[lb]{\smash{\SetFigFont{5}{6.0}{\rmdefault}{\mddefault}{\updefault}\normalsize{$c$}}}}
\put(-1274,-1861){\makebox(0,0)[lb]{\smash{\SetFigFont{5}{6.0}{\rmdefault}{\mddefault}{\updefault}\normalsize{$\sigma$}}}}
\put(12226,-5236){\makebox(0,0)[lb]{\smash{\SetFigFont{5}{6.0}{\rmdefault}{\mddefault}{\updefault}\normalsize{$a'$}}}}
\put(10426,-5236){\makebox(0,0)[lb]{\smash{\SetFigFont{5}{6.0}{\rmdefault}{\mddefault}{\updefault}\normalsize{$a$}}}}
\put(10351,-2461){\makebox(0,0)[lb]{\smash{\SetFigFont{5}{6.0}{\rmdefault}{\mddefault}{\updefault}\normalsize{$b'$}}}}
\put(12151,-2311){\makebox(0,0)[lb]{\smash{\SetFigFont{5}{6.0}{\rmdefault}{\mddefault}{\updefault}\normalsize{$e'$}}}}
\put(9601,-5011){\makebox(0,0)[lb]{\smash{\SetFigFont{5}{6.0}{\rmdefault}{\mddefault}{\updefault}\normalsize{$d$}}}}
\end{picture}

%% file: exp-contraction.pstex_t
\begin{picture}(0,0)%
\epsfig{file=exp-contraction.pstex}%
\end{picture}%
\setlength{\unitlength}{1184sp}%
\begingroup\makeatletter\ifx\SetFigFont\undefined%
\gdef\SetFigFont#1#2#3#4#5{%
  \reset@font\fontsize{#1}{#2pt}%
  \fontfamily{#3}\fontseries{#4}\fontshape{#5}%
  \selectfont}%
\fi\endgroup%
\begin{picture}(16989,7088)(-2513,-7036)
\put(901,-2761){\makebox(0,0)[lb]{\smash{\SetFigFont{5}{6.0}{\rmdefault}{\mddefault}{\updefault}\normalsize{$c'$}}}}
\put(1201,-4561){\makebox(0,0)[lb]{\smash{\SetFigFont{5}{6.0}{\rmdefault}{\mddefault}{\updefault}\normalsize{$c$}}}}
\put(1726,-5236){\makebox(0,0)[lb]{\smash{\SetFigFont{5}{6.0}{\rmdefault}{\mddefault}{\updefault}\normalsize{$a'$}}}}
\put(-74,-5236){\makebox(0,0)[lb]{\smash{\SetFigFont{5}{6.0}{\rmdefault}{\mddefault}{\updefault}\normalsize{$a$}}}}
\put(-149,-2461){\makebox(0,0)[lb]{\smash{\SetFigFont{5}{6.0}{\rmdefault}{\mddefault}{\updefault}\normalsize{$b'$}}}}
\put(1651,-2311){\makebox(0,0)[lb]{\smash{\SetFigFont{5}{6.0}{\rmdefault}{\mddefault}{\updefault}\normalsize{$e'$}}}}
\put(-899,-5011){\makebox(0,0)[lb]{\smash{\SetFigFont{5}{6.0}{\rmdefault}{\mddefault}{\updefault}\normalsize{$d$}}}}
\put(-2099,-5011){\makebox(0,0)[lb]{\smash{\SetFigFont{5}{6.0}{\rmdefault}{\mddefault}{\updefault}\normalsize{$d'$}}}}
\put(2401,-3961){\makebox(0,0)[lb]{\smash{\SetFigFont{5}{6.0}{\rmdefault}{\mddefault}{\updefault}\normalsize{$e$}}}}
\put(-749,-3961){\makebox(0,0)[lb]{\smash{\SetFigFont{5}{6.0}{\rmdefault}{\mddefault}{\updefault}\normalsize{$b$}}}}
\put(-149,-3886){\makebox(0,0)[lb]{\smash{\SetFigFont{5}{6.0}{\rmdefault}{\mddefault}{\updefault}\normalsize{$f'$}}}}
\put(1801,-3811){\makebox(0,0)[lb]{\smash{\SetFigFont{5}{6.0}{\rmdefault}{\mddefault}{\updefault}\normalsize{$f$}}}}
\put(13801,-1786){\makebox(0,0)[lb]{\smash{\SetFigFont{5}{6.0}{\rmdefault}{\mddefault}{\updefault}\normalsize{$f$}}}}
\put(13726,-7036){\makebox(0,0)[lb]{\smash{\SetFigFont{5}{6.0}{\rmdefault}{\mddefault}{\updefault}\normalsize{$a'$}}}}
\put(11926,-7036){\makebox(0,0)[lb]{\smash{\SetFigFont{5}{6.0}{\rmdefault}{\mddefault}{\updefault}\normalsize{$a$}}}}
\put(11101,-6811){\makebox(0,0)[lb]{\smash{\SetFigFont{5}{6.0}{\rmdefault}{\mddefault}{\updefault}\normalsize{$d$}}}}
\put(9901,-6811){\makebox(0,0)[lb]{\smash{\SetFigFont{5}{6.0}{\rmdefault}{\mddefault}{\updefault}\normalsize{$d'$}}}}
\put(12151,-6286){\makebox(0,0)[lb]{\smash{\SetFigFont{5}{6.0}{\rmdefault}{\mddefault}{\updefault}\normalsize{$f'$}}}}
\put(13126,-6586){\makebox(0,0)[lb]{\smash{\SetFigFont{5}{6.0}{\rmdefault}{\mddefault}{\updefault}\normalsize{$c$}}}}
\put(11026,-5686){\makebox(0,0)[lb]{\smash{\SetFigFont{5}{6.0}{\rmdefault}{\mddefault}{\updefault}\normalsize{$e'$}}}}
\put(13726,-5761){\makebox(0,0)[lb]{\smash{\SetFigFont{5}{6.0}{\rmdefault}{\mddefault}{\updefault}\normalsize{$f$}}}}
\put(14476,-5761){\makebox(0,0)[lb]{\smash{\SetFigFont{5}{6.0}{\rmdefault}{\mddefault}{\updefault}\normalsize{$e$}}}}
\put(11851,-5461){\makebox(0,0)[lb]{\smash{\SetFigFont{5}{6.0}{\rmdefault}{\mddefault}{\updefault}\normalsize{$c'$}}}}
\put(4951,-4711){\makebox(0,0)[lb]{\smash{\SetFigFont{5}{6.0}{\rmdefault}{\mddefault}{\updefault}\normalsize{Contraction}}}}
\put(4951,-1561){\makebox(0,0)[lb]{\smash{\SetFigFont{5}{6.0}{\rmdefault}{\mddefault}{\updefault}\normalsize{Deletion}}}}
\put(12901,-661){\makebox(0,0)[lb]{\smash{\SetFigFont{5}{6.0}{\rmdefault}{\mddefault}{\updefault}\normalsize{$c'$}}}}
\put(13201,-2461){\makebox(0,0)[lb]{\smash{\SetFigFont{5}{6.0}{\rmdefault}{\mddefault}{\updefault}\normalsize{$c$}}}}
\put(13726,-3136){\makebox(0,0)[lb]{\smash{\SetFigFont{5}{6.0}{\rmdefault}{\mddefault}{\updefault}\normalsize{$a'$}}}}
\put(11926,-3136){\makebox(0,0)[lb]{\smash{\SetFigFont{5}{6.0}{\rmdefault}{\mddefault}{\updefault}\normalsize{$a$}}}}
\put(13651,-211){\makebox(0,0)[lb]{\smash{\SetFigFont{5}{6.0}{\rmdefault}{\mddefault}{\updefault}\normalsize{$e'$}}}}
\put(11101,-2911){\makebox(0,0)[lb]{\smash{\SetFigFont{5}{6.0}{\rmdefault}{\mddefault}{\updefault}\normalsize{$d$}}}}
\put(9901,-2911){\makebox(0,0)[lb]{\smash{\SetFigFont{5}{6.0}{\rmdefault}{\mddefault}{\updefault}\normalsize{$d'$}}}}
\put(14401,-1861){\makebox(0,0)[lb]{\smash{\SetFigFont{5}{6.0}{\rmdefault}{\mddefault}{\updefault}\normalsize{$e$}}}}
\put(11851,-1786){\makebox(0,0)[lb]{\smash{\SetFigFont{5}{6.0}{\rmdefault}{\mddefault}{\updefault}\normalsize{$f'$}}}}
\end{picture}

%% file: exp-Tutte-triangle2.pstex_t
\begin{picture}(0,0)%
\includegraphics{exp-Tutte-triangle2.pstex}%
\end{picture}%
\setlength{\unitlength}{987sp}%
\begingroup\makeatletter\ifx\SetFigFont\undefined%
\gdef\SetFigFont#1#2#3#4#5{%
  \reset@font\fontsize{#1}{#2pt}%
  \fontfamily{#3}\fontseries{#4}\fontshape{#5}%
  \selectfont}%
\fi\endgroup%
\begin{picture}(14628,3021)(487,-5069)
\put(14551,-3361){\makebox(0,0)[lb]{\smash{{\SetFigFont{5}{6.0}{\rmdefault}{\mddefault}{\updefault}\normalsize{$c$}}}}}
\put(13651,-5011){\makebox(0,0)[lb]{\smash{{\SetFigFont{5}{6.0}{\rmdefault}{\mddefault}{\updefault}\normalsize{$a$}}}}}
\put(12751,-3361){\makebox(0,0)[lb]{\smash{{\SetFigFont{5}{6.0}{\rmdefault}{\mddefault}{\updefault}\normalsize{$b$}}}}}
\put(13651,-4411){\makebox(0,0)[lb]{\smash{{\SetFigFont{5}{6.0}{\rmdefault}{\mddefault}{\updefault}\normalsize{$\star$}}}}}
\put(8551,-3361){\makebox(0,0)[lb]{\smash{{\SetFigFont{5}{6.0}{\rmdefault}{\mddefault}{\updefault}\normalsize{$c$}}}}}
\put(7651,-5011){\makebox(0,0)[lb]{\smash{{\SetFigFont{5}{6.0}{\rmdefault}{\mddefault}{\updefault}\normalsize{$a$}}}}}
\put(7651,-4411){\makebox(0,0)[lb]{\smash{{\SetFigFont{5}{6.0}{\rmdefault}{\mddefault}{\updefault}\normalsize{$\star$}}}}}
\put(6751,-3361){\makebox(0,0)[lb]{\smash{{\SetFigFont{5}{6.0}{\rmdefault}{\mddefault}{\updefault}\normalsize{$b$}}}}}
\put(1351,-3511){\makebox(0,0)[lb]{\smash{{\SetFigFont{5}{6.0}{\rmdefault}{\mddefault}{\updefault}\normalsize{$\star$}}}}}
\put(751,-3361){\makebox(0,0)[lb]{\smash{{\SetFigFont{5}{6.0}{\rmdefault}{\mddefault}{\updefault}\normalsize{$b$}}}}}
\put(1651,-5011){\makebox(0,0)[lb]{\smash{{\SetFigFont{5}{6.0}{\rmdefault}{\mddefault}{\updefault}\normalsize{$a$}}}}}
\put(2551,-3361){\makebox(0,0)[lb]{\smash{{\SetFigFont{5}{6.0}{\rmdefault}{\mddefault}{\updefault}\normalsize{$c$}}}}}
\put(1651,-4411){\makebox(0,0)[lb]{\smash{{\SetFigFont{5}{6.0}{\rmdefault}{\mddefault}{\updefault}\normalsize{$\star$}}}}}
\end{picture}%

%% file: tour-of-tree.pstex_t
\begin{picture}(0,0)%
\epsfig{file=tour-of-tree.pstex}%
\end{picture}%
\setlength{\unitlength}{1184sp}%
\begingroup\makeatletter\ifx\SetFigFont\undefined%
\gdef\SetFigFont#1#2#3#4#5{%
  \reset@font\fontsize{#1}{#2pt}%
  \fontfamily{#3}\fontseries{#4}\fontshape{#5}%
  \selectfont}%
\fi\endgroup%
\begin{picture}(14516,4048)(-2513,-5911)
\put(-899,-5011){\makebox(0,0)[lb]{\smash{\SetFigFont{5}{6.0}{\rmdefault}{\mddefault}{\updefault}\normalsize{$d$}}}}
\put(-2099,-5011){\makebox(0,0)[lb]{\smash{\SetFigFont{5}{6.0}{\rmdefault}{\mddefault}{\updefault}\normalsize{$d'$}}}}
\put(2401,-3961){\makebox(0,0)[lb]{\smash{\SetFigFont{5}{6.0}{\rmdefault}{\mddefault}{\updefault}\normalsize{$e$}}}}
\put(-749,-3961){\makebox(0,0)[lb]{\smash{\SetFigFont{5}{6.0}{\rmdefault}{\mddefault}{\updefault}\normalsize{$b$}}}}
\put(-149,-3886){\makebox(0,0)[lb]{\smash{\SetFigFont{5}{6.0}{\rmdefault}{\mddefault}{\updefault}\normalsize{$f'$}}}}
\put(1801,-3811){\makebox(0,0)[lb]{\smash{\SetFigFont{5}{6.0}{\rmdefault}{\mddefault}{\updefault}\normalsize{$f$}}}}
\put(-149,-5236){\makebox(0,0)[lb]{\smash{\SetFigFont{5}{6.0}{\rmdefault}{\mddefault}{\updefault}\normalsize{$a$}}}}
\put(9826,-5911){\makebox(0,0)[lb]{\smash{\SetFigFont{5}{6.0}{\rmdefault}{\mddefault}{\updefault}\normalsize{\textrm{Tour of the tree}}}}}
\put(901,-2761){\makebox(0,0)[lb]{\smash{\SetFigFont{5}{6.0}{\rmdefault}{\mddefault}{\updefault}\normalsize{$c'$}}}}
\put(1201,-4561){\makebox(0,0)[lb]{\smash{\SetFigFont{5}{6.0}{\rmdefault}{\mddefault}{\updefault}\normalsize{$c$}}}}
\put(1726,-5236){\makebox(0,0)[lb]{\smash{\SetFigFont{5}{6.0}{\rmdefault}{\mddefault}{\updefault}\normalsize{$a'$}}}}
\put(-149,-2461){\makebox(0,0)[lb]{\smash{\SetFigFont{5}{6.0}{\rmdefault}{\mddefault}{\updefault}\normalsize{$b'$}}}}
\put(1651,-2311){\makebox(0,0)[lb]{\smash{\SetFigFont{5}{6.0}{\rmdefault}{\mddefault}{\updefault}\normalsize{$e'$}}}}
\end{picture}

%% file: exp-Tutte-triangle3.pstex_t
\begin{picture}(0,0)%
\includegraphics{exp-Tutte-triangle3.pstex}%
\end{picture}%
\setlength{\unitlength}{987sp}%
\begingroup\makeatletter\ifx\SetFigFont\undefined%
\gdef\SetFigFont#1#2#3#4#5{%
  \reset@font\fontsize{#1}{#2pt}%
  \fontfamily{#3}\fontseries{#4}\fontshape{#5}%
  \selectfont}%
\fi\endgroup%
\begin{picture}(14889,3021)(226,-5069)
\put(7651,-4411){\makebox(0,0)[lb]{\smash{{\SetFigFont{5}{6.0}{\rmdefault}{\mddefault}{\updefault}\normalsize{$\star$}}}}}
\put(13651,-4411){\makebox(0,0)[lb]{\smash{{\SetFigFont{5}{6.0}{\rmdefault}{\mddefault}{\updefault}\normalsize{$\star$}}}}}
\put(1651,-4411){\makebox(0,0)[lb]{\smash{{\SetFigFont{5}{6.0}{\rmdefault}{\mddefault}{\updefault}\normalsize{$\star$}}}}}
\put(6751,-5011){\makebox(0,0)[lb]{\smash{{\SetFigFont{5}{6.0}{\rmdefault}{\mddefault}{\updefault}\normalsize{$a$}}}}}
\put(8251,-2611){\makebox(0,0)[lb]{\smash{{\SetFigFont{5}{6.0}{\rmdefault}{\mddefault}{\updefault}\normalsize{$c'$}}}}}
\put(751,-5011){\makebox(0,0)[lb]{\smash{{\SetFigFont{5}{6.0}{\rmdefault}{\mddefault}{\updefault}\normalsize{$a$}}}}}
\put(2251,-2611){\makebox(0,0)[lb]{\smash{{\SetFigFont{5}{6.0}{\rmdefault}{\mddefault}{\updefault}\normalsize{$c'$}}}}}
\put(12751,-5011){\makebox(0,0)[lb]{\smash{{\SetFigFont{5}{6.0}{\rmdefault}{\mddefault}{\updefault}\normalsize{$a$}}}}}
\put(14251,-2611){\makebox(0,0)[lb]{\smash{{\SetFigFont{5}{6.0}{\rmdefault}{\mddefault}{\updefault}\normalsize{$c'$}}}}}
\put(7051,-2611){\makebox(0,0)[lb]{\smash{{\SetFigFont{5}{6.0}{\rmdefault}{\mddefault}{\updefault}\normalsize{$b$}}}}}
\put(13051,-2611){\makebox(0,0)[lb]{\smash{{\SetFigFont{5}{6.0}{\rmdefault}{\mddefault}{\updefault}\normalsize{$b$}}}}}
\put(1051,-2611){\makebox(0,0)[lb]{\smash{{\SetFigFont{5}{6.0}{\rmdefault}{\mddefault}{\updefault}\normalsize{$b$}}}}}
\put(226,-4261){\makebox(0,0)[lb]{\smash{{\SetFigFont{5}{6.0}{\rmdefault}{\mddefault}{\updefault}\normalsize{$b'$}}}}}
\put(6226,-4261){\makebox(0,0)[lb]{\smash{{\SetFigFont{5}{6.0}{\rmdefault}{\mddefault}{\updefault}\normalsize{$b'$}}}}}
\put(12226,-4261){\makebox(0,0)[lb]{\smash{{\SetFigFont{5}{6.0}{\rmdefault}{\mddefault}{\updefault}\normalsize{$b'$}}}}}
\put(3001,-4186){\makebox(0,0)[lb]{\smash{{\SetFigFont{5}{6.0}{\rmdefault}{\mddefault}{\updefault}\normalsize{$c$}}}}}
\put(9001,-4186){\makebox(0,0)[lb]{\smash{{\SetFigFont{5}{6.0}{\rmdefault}{\mddefault}{\updefault}\normalsize{$c$}}}}}
\put(15001,-4186){\makebox(0,0)[lb]{\smash{{\SetFigFont{5}{6.0}{\rmdefault}{\mddefault}{\updefault}\normalsize{$c$}}}}}
\put(8026,-3511){\makebox(0,0)[lb]{\smash{{\SetFigFont{5}{6.0}{\rmdefault}{\mddefault}{\updefault}\normalsize{$\star$}}}}}
\put(8326,-5011){\makebox(0,0)[lb]{\smash{{\SetFigFont{5}{6.0}{\rmdefault}{\mddefault}{\updefault}\normalsize{$a'$}}}}}
\put(14326,-5011){\makebox(0,0)[lb]{\smash{{\SetFigFont{5}{6.0}{\rmdefault}{\mddefault}{\updefault}\normalsize{$a'$}}}}}
\put(2326,-5011){\makebox(0,0)[lb]{\smash{{\SetFigFont{5}{6.0}{\rmdefault}{\mddefault}{\updefault}\normalsize{$a'$}}}}}
\end{picture}%

%% file: case-study-external2.pstex_t
\begin{picture}(0,0)%
\epsfig{file=case-study-external2.pstex}%
\end{picture}%
\setlength{\unitlength}{1579sp}%
\begingroup\makeatletter\ifx\SetFigFont\undefined%
\gdef\SetFigFont#1#2#3#4#5{%
  \reset@font\fontsize{#1}{#2pt}%
  \fontfamily{#3}\fontseries{#4}\fontshape{#5}%
  \selectfont}%
\fi\endgroup%
\begin{picture}(9150,2813)(8851,-3211)
\put(17176,-3211){\makebox(0,0)[lb]{\smash{\SetFigFont{5}{6.0}{\rmdefault}{\mddefault}{\updefault}\normalsize{$v_0$}}}}
\put(17851,-961){\makebox(0,0)[lb]{\smash{\SetFigFont{5}{6.0}{\rmdefault}{\mddefault}{\updefault}\normalsize{$h'$}}}}
\put(17851,-2461){\makebox(0,0)[lb]{\smash{\SetFigFont{5}{6.0}{\rmdefault}{\mddefault}{\updefault}\normalsize{$h$}}}}
\put(16501,-961){\makebox(0,0)[lb]{\smash{\SetFigFont{5}{6.0}{\rmdefault}{\mddefault}{\updefault}\normalsize{$h_*'$}}}}
\put(16576,-2461){\makebox(0,0)[lb]{\smash{\SetFigFont{5}{6.0}{\rmdefault}{\mddefault}{\updefault}\normalsize{$h_*$}}}}
\put(16426,-1711){\makebox(0,0)[lb]{\smash{\SetFigFont{5}{6.0}{\rmdefault}{\mddefault}{\updefault}\normalsize{$e_*$}}}}
\put(18001,-1711){\makebox(0,0)[lb]{\smash{\SetFigFont{5}{6.0}{\rmdefault}{\mddefault}{\updefault}\normalsize{$e$}}}}
\put(8926,-961){\makebox(0,0)[lb]{\smash{\SetFigFont{5}{6.0}{\rmdefault}{\mddefault}{\updefault}\normalsize{$h_*'$}}}}
\put(8851,-1711){\makebox(0,0)[lb]{\smash{\SetFigFont{5}{6.0}{\rmdefault}{\mddefault}{\updefault}\normalsize{$e_*$}}}}
\put(9001,-2461){\makebox(0,0)[lb]{\smash{\SetFigFont{5}{6.0}{\rmdefault}{\mddefault}{\updefault}\normalsize{$h_*$}}}}
\put(10426,-1711){\makebox(0,0)[lb]{\smash{\SetFigFont{5}{6.0}{\rmdefault}{\mddefault}{\updefault}\normalsize{$e$}}}}
\put(9601,-3211){\makebox(0,0)[lb]{\smash{\SetFigFont{5}{6.0}{\rmdefault}{\mddefault}{\updefault}\normalsize{$v_0$}}}}
\put(10276,-961){\makebox(0,0)[lb]{\smash{\SetFigFont{5}{6.0}{\rmdefault}{\mddefault}{\updefault}\normalsize{$h'$}}}}
\put(10276,-2461){\makebox(0,0)[lb]{\smash{\SetFigFont{5}{6.0}{\rmdefault}{\mddefault}{\updefault}\normalsize{$h$}}}}
\end{picture}